\documentclass[11pt]{article}
\usepackage[T2A]{fontenc}
\usepackage[cp1251]{inputenc}
\usepackage{amsfonts}
\usepackage{eufrak}
\usepackage{amssymb}
\usepackage{amsmath}

\begin{document}

\noindent{\large\textbf{{On the Skitovich--Darmois theorem for some
}}}

\noindent{\large\textbf{{locally compact Abelian groups
}}}

\bigskip

\noindent{Gennadiy Feldman and Margaryta Myronyuk}

\bigskip

\noindent{\bf Abstract.}  Let
$X$ be a locally compact Abelian group, $\alpha_{j}, \beta_j$ be  topological automorphisms of $X$. Let $\xi_1, \xi_2$ be independent random variables with
values in $X$ and distributions $\mu_j$ with non-vanishing characteristic functions.
It is known that if $X$ contains no subgroup topologically isomorphic to the circle group
$\mathbb{T}$, then the independence of the linear forms  $L_1=\alpha_1\xi_1+\alpha_2\xi_2$
and $L_2=\beta_1\xi_1+\beta_2\xi_2$ implies that  $\mu_j$ are Gaussian distributions. We  prove  that
if  $X$ contains no subgroup topologically isomorphic to
$\mathbb{T}^2$, then the independence of   $L_1$ and $L_2$ implies that
$\mu_j$ are either Gaussian distributions
or convolutions of Gaussian distributions and signed measures
 supported in a subgroup of  $X$ generated by an element of order  2. The proof is based on solving the Skitovich--Darmois functional equation on some locally compact Abelian groups.

\bigskip


\bigskip


\bigskip

\noindent{\textbf{1. Introduction}}

\bigskip

\noindent One of the most well-known characterization  theorems in mathematical statistics is the following
theorem which characterizes the Gaussian distribution
on the real line.

\textbf{The Skitovich--Darmois theorem} (\cite[Ch. 3]{Kag-Lin-Rao}).
\textit{Let $\xi_j, \ j=1, 2,\dots, n,\ n\geq 2$,
be independent random variables, $\alpha_j, \beta_j$ be nonzero real numbers.
If the linear forms   $L_1=\alpha_1\xi_1+\cdots+\alpha_n\xi_n$ and
$L_2=\beta_1\xi_1+\cdots+\beta_n\xi_n$ are independent, then all
 $\xi_j$ are Gaussian.}

The Skitovich--Darmois theorem was generalized by S.G. Ghurye and I. Olkin
 to the case when instead of random variables random vectors $\xi_j$  in the space $\mathbb{R}^a$ are considered
and coefficients of the linear forms      $L_1$ and $L_2$ are nonsingular matrices. They proved that in this
case the independence of
$L_1$ and $L_2$ implies that all   $\xi_j$ are also Gaussian (\cite[Ch. 3]{Kag-Lin-Rao}). The Skitovich--Darmois theorem was generalized in different directions (see e.g. \cite{NeueScho}, \cite{ChisGets}, \cite{My-2008}, \cite{Ibrahim}).
 Especially many publications have been devoted to
group analogues of the Skitovich--Darmois theorem in the case, when independent random variables take values in a
locally compact Abelian group, and coefficients of the forms are topological
automorphisms of the group (see e.g. \cite{SMZh}, \cite{Fe-SD-2003},   \cite{JOTP2015}--\cite{MiFe5}, \cite{Ma}, \cite{MiFe4},  and also \cite{Fe1}, \cite{Fe5}  where one can find
additional  references). In this paper  we continue these research.
It should be noted that the study of group analogues of the Skitovich--Darmois theorem on a locally compact Abelian group $X$ is based on the study of
solutions of the Skitovich--Darmois functional equation on the character group of
$X$.

Denote by   $\mathbb{T}$ the circle group (the one dimensional torus), i.e.
${\mathbb{T}=\{z\in\mathbb{C}:\ |z|=1\}}$. The following theorem was proved in \cite{Fe-SD-2003}.

\textbf{Theorem A.} \textit{ Let $X$ be a second countable locally compact Abelian group containing no subgroups
topologically isomorphic to  $\mathbb{T}$. Let $\alpha_j$, $\beta_j$,  $j=1,2,\dots,n$, $n\geq 2$, be
topological automorphisms of the group $X$. Let $\xi_j$ be
independent random variables with values in  $X$ and
 distributions  $\mu_j$ with non-vanishing characteristic
functions. Then the independence of the linear
forms  $L_1 = \alpha_1\xi_1 + \cdots + \alpha_n\xi_n$ and $L_2 =
\beta_1\xi_1 + \cdots + \beta_n\xi_n$ implies that all
 $\mu_j$ are Gaussian distributions}.

As  noted for example in \cite{Fe-Be-1986}, if a locally compact Abelian group $X$  contains a subgroup topologically isomorphic to
$\mathbb{T}$, Theorem  A   fails even for the simplest linear forms
$L_1=\xi_1+\xi_2$ and  $L_2=\xi_1-\xi_2$. Thus Theorem A
gives a complete description of locally compact Abelian groups for which the Skitovich--Darmois
theorem holds in its classical formulation under assumption that the characteristic functions
of  $\mu_j$  do not vanish. We will formulate now the following general problem.

\textbf{Problem 1.} \textit{ Let $X$ be a second countable locally compact Abelian group, $\alpha_j$, $\beta_j$, $j = 1, 2,\dots, n$, $n \ge 2,$ be
topological automorphisms of  $X$.
Let $\xi_j$ be
independent random variables with values in the group $X$ and
 distributions  $\mu_j$ with non-vanishing characteristic
functions. Assume that  the linear
forms  $L_1 = \alpha_1\xi_1 + \cdots + \alpha_n\xi_n$ and $L_2 =
\beta_1\xi_1 + \cdots + \beta_n\xi_n$ are independent. Describe the possible distributions
 $\mu_j$.}

Taking into account Theorem  A, it is sufficient to solve Problem 1 for groups $X$ containing a subgroup
topologically isomorphic to   $\mathbb{T}$. Problem 1  for arbitrary $n\geq 2$ is complicated enough. For  $n=2$
a partial solution of Problem   1 for the group  $X=\mathbb{T}^2$
was obtained in  \cite{MiFe4}, and a complete solution for the groups
 $X=\mathbb{R}\times \mathbb{T}$ and
$X=\Sigma_\text{\boldmath $a$}\times \mathbb{T}$, where
$\Sigma_{\text{\boldmath $a$}}$ is an  $\text{\boldmath $a$}$-adic solenoid, was obtained in  \cite{MiFe5}.

The main result of the article is a complete solution
of Problem 1 for  $n=2$
for an arbitrary locally compact Abelian group  $X$, containing a  subgroup
topologically isomorphic to  $\mathbb{T}$, but not containing a subgroup
topologically isomorphic to $\mathbb{T}^2$.
We prove that in this case the distributions  $\mu_j$
are not only Gaussian, but convolutions of Gaussian distributions
and some signed measures  (Theorem 3). The proof is based on the fact that
solution of Problem 1 for an arbitrary locally compact Abelian group
can be reduced to the case when a group is of the form
$\mathbb{R}^{\mathfrak{n}} \times \mathbb{T}^{\mathfrak{m}}$, where ${\mathfrak{n}}\le
\aleph_0$,   ${\mathfrak{m}}\le \aleph_0$ (Theorems 1 and 2).
As to the case when a locally compact Abelian group  $X$  contains a  subgroup
topologically isomorphic to $\mathbb{T}^2$, taking into account results of
\cite{MiFe4}, we can hardly hope for a complete solution of Problem 1.

We will use in the article the standard results of abstract harmonic analysis (see \cite{Hewitt-Ross}). Remind some definitions and agree on notation.
Let $X$ be a locally compact Abelian group,
   $Y=X^*$ be its character group,
$(x,y)$ be the value of a character  $y\in Y$ at an element  $x\in X$,
${\rm Aut}(X)$ be the group of topological automorphisms of
$X$. A Borel subgroup of the group  $X$ is called characteristic if it is invariant with respect
to any automorphism $\delta\in{\rm Aut}(X)$.
Let $A$ be a closed subgroup of $X$, and  $\delta$ be a continuous endomorphism of $X$. Denote by
$\delta|_A$ the restriction of $\delta$ to  $A$.
If   $H$ is a subgroup of the group $Y$, denote by   $A(X,H)=\{x \in X:
(x,y)=1$ for all  $y \in H\}$ its annihilator. Let $X_1$ and  $X_2$ be  locally compact Abelian groups with character groups
$Y_1$ and  $Y_2$ respectively. For any continuous homomorphism
 $f : X_1 \mapsto X_2$  define an adjoint homomorphism
 $\tilde{f} :Y_2 \mapsto Y_1$ by the formula
$(x_1,\tilde{f}y_2)=(f x_1,y_2)$ for all  $x_1 \in X_1$, $y_2
\in Y_2$. Denote by $c_X$  the connected component of zero of
the group $X$. Denote by $\mathbb{Z}$ the group of integers, by
 $\mathbb{Z}(2)$ the subgroup of $\mathbb{T}$ consisting of elements
 $\pm 1$. Denote by  $\mathbb{T}^{{\mathfrak{m}}}$, where
${{\mathfrak{m}}}\le {\aleph_0}$,   the finite-dimensional  or infinite-dimensional torus. Denote by
  $\mathbb{Z}^{{\mathfrak{m}}*}$ the free group of rank
   ${\mathfrak{m}}$.  Denote by $I$ the identity automorphism of a group. Put $Y^{(2)}=\{2y: \ y\in
Y\}$.

 Denote by   $\mathbb{R}^{\aleph_0}$ the space of all sequences of real numbers in the topology
 of  the projective limit of  spaces  $\mathbb{R}^{n}$ (in the product topology), and by
    $\mathbb{R}^{\aleph_0*}$ the space of all finite sequences
 of real numbers with the topology of  strictly inductive limit of spaces
    $\mathbb{R}^{n}$. We note that the groups $\mathbb{R}^{\aleph_0}$ and
$\mathbb{R}^{\aleph_0*}$ are the character group of one another.

Let  $\psi(y)$ be a function on $Y$, and $h$ be an arbitrary element of $Y$. Denote by  $\Delta_h$
the finite difference operator
$$\Delta_h \psi(y)=\psi(y+h)-\psi(y), \quad y\in Y.$$
Denote by  $\langle . , . \rangle$  the scalar product in
the space $\mathbb{R}^{b}$.

Let ${M^1}(X)$ be the convolution semigroup of probability distributions on
   $X$ and $\mu \in {M^1}(X)$.  Denote by
 $$\hat \mu(y) = \int_X (x, y) d\mu(x)$$
the characteristic function (Fourier transform) of $\mu$, and by
 $\sigma(\mu)$ the support of  $\mu$.
For $\mu \in {M}^1(X)$ define the distribution  $\bar \mu \in
{M}^1(X)$ by the formula   $\bar \mu(B) = \mu(-B)$
for any Borel set   $B$. We note that
$\hat{\bar{\mu}}(y)=\overline{\hat\mu(y)}$.
Denote by  $E_x$ the degenerate distribution concentrated at a point $x
\in X$.

A distribution $\gamma$ on the group  $X$ is called Gaussian in the sense of Parthasarathy
  (\cite[Ch. IV]{Parthasarathy},  if its characteristic function is of the form
$$
\hat\gamma(y)=(x,y)\exp\{-\varphi (y)\},
$$
where $x\in X$, and  $\varphi(y)$
is a continuous nonnegative function on     $Y$, satisfying the equation
\begin{equation}
\label{1} \varphi (u+v)+\varphi (u-v)=2[\varphi (u)+\varphi
(v)],\quad u,v\in Y.
\end{equation}
Since  we will deal only with  Gaussian distributions in the sense of Parthasarathy
we will call them Gaussian. Denote by $\Gamma(X)$  the set of Gaussian distributions on
the group $X$.

\bigskip

\noindent\textbf{2. Reduction of Problem 1 for arbitrary groups to the groups of the form
    $\mathbb{R}^{\mathfrak{n}} \times \mathbb{T}^{\mathfrak{m}}$}

\bigskip

 It is convenient for us to formulate as lemma the following statement.

\textbf{Lemma 1} (\cite[\S10.1]{Fe1}). {\it  Let $X$ be a second countable locally compact Abelian group, $\alpha_j, \ \beta_j
\in {\rm Aut}(X)$.
Let $\xi_j, \ j = 1, 2,\dots, n, \ n \ge 2,$ be
independent random variables with values in  $X$ and
 distributions  $\mu_j$. The  linear forms
$L_1=\alpha_1\xi_1+\cdots+\alpha_n\xi_n$
$L_2=\beta_1\xi_1+\cdots+\beta_n\xi_n$  are independent if and only
if the characteristic functions $\hat\mu_j(y)$  satisfy the Skitovich--Darmois functional equation}
\begin{equation}
\label{2} \prod_{j = 1}^{n}{\hat\mu_j(\tilde\alpha_j u +
\tilde\beta_j v)} = \prod_{j = 1}^{n}
{\hat\mu_j(\tilde\alpha_j u)} \prod_{j = 1}^{n}
{\hat\mu_j(\tilde\beta_j v)}, \quad u,v \in Y.
\end{equation}

Taking into account that the characteristic function of the distribution $\mu_j$
is the mathematical expectation $\hat\mu_j(y)=${\bf
E}$[(\xi_j,y)]$, the proof of Lemma  1 is the same as in the classical case
    $X=\mathbb{R}$.

We need the following K. Stein theorem which we formulate as a lemma.

\textbf{Lemma 2} (\cite[\S19.3]{Fu1}).\textit{ Let $H$ be a countable Abelian group. Then
   $H$ is represented in the form
   $H=N\times M$, where $M$ is  free, and
 $N$ has no free factor-groups. The group $N$ is uniquely determined by the group $H$}.

Lemma  2 implies the following statement.

\textbf{Lemma 3.} {\it  Let $X$ be a second countable connected locally compact Abelian group.
Then the group  $X$ is topologically isomorphic to a group of the form  $\mathbb{R}^a\times
K\times \mathbb{T}^{\mathfrak{m}}$, where $a\geq 0$, $K$ is a connected compact Abelian group containing
no subgroup topologically isomorphic to
$\mathbb{T}$, ${\mathfrak{m}}\le \aleph_0$.  Furthermore, $a$ and  ${\mathfrak{m}}$
are uniquely determined by the group  $X$,  and the group  $K$ is uniquely determined by the group  $X$ up to a topological isomorphism.}

\textbf{Proof.} By the structure theorem for connected locally compact Abelian groups
    $X=L\times F$,  where
 $L\cong\mathbb{R}^a,$
$a\ge 0$, and  $F$ is a second countable connected compact Abelian group. Furthermore,
  $a$ and $F$ are uniquely determined by the group   $X$.
Put  $H=F^*$. Then $H$ is
a countable discrete  torsion free Abelian group. By Lemma  2  $H$  is represented in the form
$H=N\times M$, where $M$ is a free Abelian group,   $N$ is a  torsion free Abelian
group without free factor-groups.
Furthermore, the group    $N$ is uniquely determined by the group  $H$. This implies that $F=K\times G$,
where $K\cong N^*$ is a connected compact group containing no subgroup topologically isomorphic to  $\mathbb{T}$, $G\cong M^*\cong
\mathbb{T}^{\mathfrak{m}}$,  where ${\mathfrak{m}}\le \aleph_0$.
Since $G=A(F, N)$,  the group $G$ is uniquely determined by the group   $F$, and hence by the group
$X$ too. It follows from this that
${\mathfrak{m}}$ is uniquely determined by the group  $X$, and the subgroup $K$, which is topologically
isomorphic to the factor-group $F/G$,  is uniquely determined up to a topological isomorphism.  $\Box$

We also need the following standard lemma.

 \textbf{Lemma 4} (see e.g. \cite[\S2.5]{Fe}).\textit{ Let $X$ be
a topological Abelian group, $B$ be a Borel subgroup of $X$, $\mu$ be a distribution on   $X$ concentrated
on   $B$. Let
$\mu=\mu_1*\mu_2$, where  $\mu_j$ are distributions on   $X$.
Then the distributions $\mu_j$ can be replaced by their
shifts $\mu'_j$ in such a manner that $\mu=\mu'_1*\mu'_2$ and $\mu'_j$ are  concentrated
on $B$.}

\textbf{Lemma 5} (\cite{Fe-SD-2003}).\textit{ Let $X$ be a second countable locally compact Abelian group, $\alpha_j, \ \beta_j
\in {\rm Aut}(X)$.
Let $\xi_j$, $j=1,2,\dots,n$, $n\geq
2$,  be
independent random variables with values in  $X$ and
 distributions  $\mu_j$ with non-vanishing characteristic
functions. Then the independence of the linear
forms  $L_1 = \alpha_1\xi_1 + \cdots + \alpha_n\xi_n$ and $L_2 =
\beta_1\xi_1 + \cdots + \beta_n\xi_n$ implies that there exist  elements
$x_j\in X$,   $j=1, 2, \dots, n$ such that all distributions
${\mu}'_j$ of random variables $\xi'_j=\xi_j+x_j$ are supported in
 the subgroup $c_X$.}

Taking into account that   $c_X$ is a characteristic subgroup,
it follows from Lemma 5 that
the study of distributions of independent random variables  with non-vanishing characteristic
functions which are characterized by the independence of the linear forms
     $L_1$ and
$L_2$ is reduced to the case when $X$ is a connected group. Although
the structure of connected locally compact Abelian groups is still complicated
  (each such group is topologically isomorphic to a group of the form
  $\mathbb{R}^a\times F,$   $a\ge 0$,  $F$
is a connected compact Abelian group), it is much simpler
than the structure of
arbitrary locally compact Abelian groups.

In this section we considerably strengthen Lemma 5.  Namely,
we prove that      the study of distributions of independent
random variables  with non-vanishing characteristic
functions which are characterized by the independence of linear forms
     $L_1$ and
$L_2$ is reduced to the case when
the group   $X$  is topologically isomorphic to a group of the form
$\mathbb{R}^{\mathfrak{n}} \times \mathbb{T}^{\mathfrak{m}}$,
where ${\mathfrak{n}}\le \aleph_0$,  ${\mathfrak{m}}\le \aleph_0$.
It should be noted that if   ${\mathfrak{n}}=\aleph_0$,  then
$X$ is not a locally compact group, although its structure is
simple enough.

 Let $X$ be a second countable  locally compact Abelian group. By Lemma 3 $c_X\cong\mathbb{R}^a\times
K\times \mathbb{T}^{\mathfrak{m}}$, where $a\geq 0$, $K$ is a connected compact Abelian group containing no subgroup
topologically isomorphic to $\mathbb{T}$, ${\mathfrak{m}}\le \aleph_0$.
Since the group
 $K$ is uniquely determined by the group $X$ up to a topological
 isomorphism, the dimension of the group
   $K$ is uniquely determined by the group $X$. First consider
   the case when  $K$
has a finite dimension.

\textbf{Theorem 1.} \textit{Let $X$ be a second countable  locally compact Abelian group such
that its connected component of zero  $c_X\cong\mathbb{R}^a\times K\times
\mathbb{T}^{\mathfrak{m}}$, where $a\geq 0$, $K$ be
a connected compact Abelian group containing no subgroup
topologically isomorphic to
 $\mathbb{T}$, ${\mathfrak{m}}\le \aleph_0$.
Assume that   $K$ has a finite dimension   $l$.
Let $\alpha_j, \ \beta_j \in \mathrm{Aut}(X)$.
Let $\xi_j$, $j = 1, 2,\dots, n$, $n \ge 2,$ be
independent random variables with values in  $X$ and
 distributions  $\mu_j$ with non-vanishing characteristic
functions.
Assume that the linear forms   $L_1 = \alpha_1\xi_1 + \cdots +
\alpha_n\xi_n$  and $L_2 = \beta_1\xi_1 + \cdots + \beta_n\xi_n$
are independent.   Then there exist  a continuous monomorphism
$p:G\mapsto X$,  where  $G=
\mathbb{R}^{b}\times\mathbb{T}^{\mathfrak{m}}$, $b=a+l$, and elements
 $x_j\in X$, \ $j=1, 2, \dots, n,$ such that all distributions
 $\mu_j*E_{x_j}$ are concentrated on the subgroup
  $p(G)$. Furthermore, $p(G)$ is a characteristic subgroup.}

\textbf{Proof.}  It follows from Lemma 5 that we can assume from
the beginning that $X$ is a connected group,
 i.e.  $X=c_X$. Then $Y\cong \mathbb{R}^a\times D\times \mathbb{Z}^{\mathfrak{m}*}$, where $D=K^*$
is a countable discrete torsion free Abelian group containing no free factor-groups. To avoid introducing
additional notation we assume that
$X=\mathbb{R}^a\times K\times \mathbb{T}^{\mathfrak{m}}$ and
$Y=\mathbb{R}^a\times D\times \mathbb{Z}^{{\mathfrak{m}}*}$.
Since the group  $K$ has dimension  $l$,  the rank of the group  $D$
is also $l$. Put $G=\mathbb{R}^{b}\times
\mathbb{T}^{\mathfrak{m}}$,  $b=a+l$. Then  $H=G^*\cong\mathbb{R}^{b}\times
\mathbb{Z}^{{\mathfrak{m}}*}$. We also assume that $H=\mathbb{R}^{b}\times
\mathbb{Z}^{{\mathfrak{m}}*}$. Denote by
$(s,d,k)$,
$s\in \mathbb{R}^a$, $d\in D$, $k\in
\mathbb{Z}^{{\mathfrak{m}}*}$, elements of the group  $Y$, and by $(w, k)$,   $w\in \mathbb{R}^{b}$, $k\in
\mathbb{Z}^{{\mathfrak{m}}*}$, elements of the group  $H$.
Construct a natural embedding of the group $Y$ to the group  $H$ (see e.g. \cite[\S5.6]{Fe}).
Choose in $D$ a maximal independent
system of elements  $d_1, \dots, d_l$.  Then for every $d \in D$
there exist integers
  $q\neq 0, q_1, \dots, q_l$, such that
\begin{equation}
\label{abc1}qd=q_1 d_1+\cdots +q_l d_l.
\end{equation}
The independence of the set $\{d_1, \dots, d_l\}$ implies that the
rational numbers
   $\{q_j/q\}$   are uniquely determined by
 $d$.
Define the mapping   $f_0:{\mathbb{R}^a\times D\mapsto
\mathbb{R}^b}$ by the formula $f_0(s, d)=(s, q_1/q,\dots,q_l/q)$, $s\in
\mathbb{R}^a,$ $d\in D.$  Since  $D$ is a torsion free group,
 $f_0$   is a continuous monomorphism of the group
${\mathbb{R}^a\times D}$ to $\mathbb{R}^b$.
Extend  $f_0$ from  ${\mathbb{R}^a\times D}$
 to the continuous monomorphism   $f:Y\mapsto H$  putting
\begin{equation}
\label{abc2}f y=f(s,d,k)=(s, q_1/q,\dots,q_l/q,k), \quad
y=(s,d,k), \quad s\in \mathbb{R}^a, \ d\in D, \ k\in
\mathbb{Z}^{{\mathfrak{m}}*}.
\end{equation}
Since the group   $K$ has no subgroups topologically isomorphic to
$\mathbb{T}$, we have $\overline{f_0({\mathbb{R}^a\times
D})}=\mathbb{R}^b$ (\cite{TVP1979}). Then obviously
$\overline{f(Y)}=H$.  Set $p=\tilde{f}$. Since
$\overline{f(Y)}=H$, it follows from the properties of adjoint homomorphisms that
 $p:G\mapsto X$ is a continuous monomorphism.

Note now that the inequality
\begin{equation}
\label{11}|\psi(u)-\psi(v)|^2\le2(1-{\rm Re}\
\psi(u-v)), \quad u, v\in Y,
\end{equation}
holds true for any positive definite function
$\psi(y)$    on an arbitrary Abelian group
 $Y$.

By Lemma  1 the characteristic functions  $\hat\mu_j(y)$  satisfy
the Skitovich--Darmois functional equation (\ref{2}).
Consider the distributions   $\nu_j=\mu_j*\bar\mu_j, \  \ j=1,2,\dots,n$.
Then $\hat\nu_j(y)=|\hat\mu_j(y)|^2> 0$. Obviously, the characteristic functions
  $\hat\nu_j(y)$  also satisfy the Skitovich--Darmois functional equation
 (\ref{2}).

Since  $\mathbb{R}^a$ is a connected component of zero of the group     $Y$,
it follows that $\mathbb{R}^a$ is a characteristic subgroup of the group    $Y$.
For this reason  $K\times \mathbb{T}^{\mathfrak{m}}=
A(X,\mathbb{R}^a)$ is
a characteristic subgroup of the group
   $X$. It follows from Lemma  2 that
$D$  is a characteristic subgroup of the group $D\times
\mathbb{Z}^{{\mathfrak{m}}*}$. This implies that
$\mathbb{T}^{\mathfrak{m}}=A(K\times \mathbb{T}^{\mathfrak{m}},
D)$  is a characteristic subgroup of the group    $ K\times
\mathbb{T}^{\mathfrak{m}}$, and hence, $\mathbb{T}^{\mathfrak{m}}$ is a characteristic subgroup of the group $X$.
It follows from this that    $\mathbb{R}^a\times D=A(Y,
\mathbb{T}^{\mathfrak{m}})$   is a characteristic subgroup of the group
 $Y$, and we can consider the restriction of the Skitovich--Darmois functional equation (\ref{2})
for the characteristic functions     $\hat\nu_j(y)$
to this subgroup. Since  $\mathbb{R}^a\times D\cong (\mathbb{R}^a\times
K)^*$,  and the group $\mathbb{R}^a\times K$ has no subgroups
topologically isomorphic to
  $\mathbb{T}$, it follows from Lemma   1 and Theorem  A that these restrictions
  are the characteristic functions of some Gaussian distributions
  on the group $\mathbb{R}^a\times K$. Taking into account that
$\hat\nu_j(y)>0$, we have the representations
    \begin{equation}
\label{new1}\hat\nu_j(y)=\exp\{-\varphi_j (y)\}, \quad  y\in
\mathbb{R}^a\times D, \quad    j=1,2,\dots,n,
\end{equation}
where  $\varphi_j(y)$ are continuous nonnegative functions on
$\mathbb{R}^a\times D$ satisfying equation   (\ref{1}).
As has been proved in  \cite[\S5.6]{Fe}, it follows
from the properties of the functions   $\varphi_j(y)$
that there exist symmetric positive semidefinite
 $(b\times b)$-matrices $Q_j$, such that
\begin{equation}
\label{new2}\varphi_j (y)=\langle Q_j f_0 y,f_0 y\rangle, \quad
y\in \mathbb{R}^a\times D,   \quad    j=1,2,\dots,n.
\end{equation}
Assume that   $(w, k)\in f(Y)$, i.e.  $(w, k)=f(y)$, $y\in
Y$. Consider on the subgroup   $f(Y)$ the functions $h_j(w,
k)=\hat\nu_j(f^{-1}(w, k))$,  $j=1, 2,\dots, n.$ Since
$f(y)=f_0(y)$ for $y\in  \mathbb{R}^a\times D$, it follows
from  (\ref{new1})
and (\ref{new2}) that
\begin{equation}
\label{1e1}h_j(w, 0)=\exp\{-\langle Q_j
w, w\rangle\},    \quad  (w, 0)\in f(\mathbb{R}^a\times D), \quad  j=1,2,\dots,n.
\end{equation}

Taking into account (\ref{11}) it follows from (\ref{1e1}) that positive definite functions
 $h_j(w, k)$  are   uniformly continuous on the subgroup   $f(Y)$
in topology induced on $f(Y)$  by the standard topology of
$H$. Since $\overline{f(Y)}=H$, the functions
   $h_j(w, k)$ can be extended  by continuity to some continuous functions
     $\bar h_j(w, k)$,
$w\in \mathbb{R}^{b}$, $k\in \mathbb{Z}^{{\mathfrak{m}}*}$,
on the group $H$. Obviously,  $\bar h_j(w, k)$ are also positive definite functions.
By the Bochner theorem there exist distributions $\lambda_j$ on   $G$ such that
$\hat\lambda_j(w, k)=\bar h_j(w, k)$,   $w\in \mathbb{R}^{b}$,
$k\in \mathbb{Z}^{{\mathfrak{m}}*}$. Since
$\hat\lambda_j(f(y))=\hat\nu_j(y)$ for all  $y\in Y$,
then $\nu_j=p(\lambda_j)$. Hence, the distributions  $\nu_j$
are concentrated on $p(G)$. It is obvious that $p(G)$ is a Borel
subgroup.
By Lemma 4 this implies that there exist elements
 $x_j\in X$, \ $j=1, 2, \dots, n,$ such that all distributions
 $\mu_j*E_{x_j}$  are concentrated on  $p(G)$.

It remains to prove that   $p(G)$ is a characteristic subgroup.
Let $\delta\in{\rm Aut}(X)$. Put
$A_{{\delta}}=f\tilde{\delta}f^{-1}$. Then $A_{{\delta}}$ is an algebraic
 automorphism of the subgroup $f(Y)$. The automorphism
${\tilde{\delta}}$   is determined by its restriction on the subgroup
$\mathbb{R}^a$ and its values on the maximal independent system elements of the subgroup
   $ D\times \mathbb{Z}^{\mathfrak{m}*}$.
For this reason the automorphism  ${\tilde{\delta}}$ determines in a natural way a matrix $A(\delta)$.
 It is easily seen that  if $f y=(s, q_1/q,\dots,q_l/q,k), \ y=(s,d,k),
\ s\in \mathbb{R}^a, \ d\in D, \ k\in
\mathbb{Z}^{{\mathfrak{m}}*}$, then  $A_\delta f y=A(\delta)fy$,
where the expression in the right-hand side of this equality is
the product of the matrix  $A(\delta)$
  and the vector  $fy$. It is clear that the automorphism
$A_{{\delta}}$   can be uniquely extended  to the topological
automorphism of
the group   $H$.  Denote by
   $\bar A_{{\delta}}$ this extended automorphism. We have
$\tilde{\delta} y = f^{-1}{A_{{\delta}}} f y$, ${A_{{\delta}}}
f y={\bar A_{{\delta}}} f  y, \ y\in Y$. Let $g\in G,$ $y\in Y$.
Then  $(\delta p g, y)=(p g, \tilde{\delta} y)= (p
g,f^{-1}A_{{\delta}} f y)= (g, A_{{\delta}} f y)=(g, \bar
A_{{\delta}} f y)=(p \widetilde{\bar{A}_{{\delta}}} g, y)$. It
follows from this that
  $\delta p g=p \widetilde{\bar{A}_{{\delta}}} g$.
Hence $p(G)$ is a characteristic subgroup.  $\Box$

\textbf{ Remark 1.} We keep the notation used in the proof of
Theorem 1. Let
   $X=\mathbb{R}^a\times K\times
\mathbb{T}^{\mathfrak{m}}$. Denote by   $L$
the arcwise connected component of zero of the group $X$. We verify that
$p(G)=L$. By  the Dixmier theorem   $L$ is a union of all one-parametric subgroups of the group
     $X$
(\cite[\S 8.19]{Ar}).  Let $p_1:\mathbb{R}\mapsto X$ be a continuous homomorphism. Put $f_1=\tilde p_1$.
Taking into account  Dixmier's theorem the required statement
will be proved if we check that that
 $p_1 t \in p(G)$  for all $t\in \mathbb{R}$.  To this end
 it suffices to show that there exists an element
  $g_t\in G$ such that
$(y, p_1 t )=(y, p g_t )$ for all  $y\in Y$. Let
$\{e_j\}_{j=1}^a$ be a natural basis in   $\mathbb{R}^a$, and
$\{b_j\}$  be a natural basis in    $\mathbb{Z}^{\mathfrak{m}*}.$
Taking into account (\ref{abc1}), we obtain
\begin{equation}
\label{abc3}f_1 y =f_1(s,d,k)=\sum_{j=1}^as_jf_1 e_j
+\sum_{j=1}^l{q_j\over q}f_1 d_j +\sum_{j}k_jf_1 b_j,
$$
$$
s=(s_1, \dots, s_a)\in \mathbb{R}^a, \ d\in D, \ k=(k_1, k_2,
\dots)\in \mathbb{Z}^{{\mathfrak{m}}*}.$$

This implies that

 $$(y, p_1 t )=
(f_1 y , t)=\exp\left\{it\left(\sum_{j=1}^as_jf_1 e_j
+\sum_{j=1}^l{q_j\over q}f_1 d_j +\sum_{j}k_jf_1 b_j
\right)\right\}.
\end{equation}
Put  $$g_t=(tf_1e_1,\dots, tf_1e_a, tf_1d_1,\dots, tf_1d_l,
e^{itf_1b_1},e^{itf_1b_2},\dots)\in G.$$ Let $y=(s_1, \dots, s_a,
d, k_1, k_2, \dots)\in Y.$ Taking into account (\ref{abc2}) we get
\begin{equation}
\label{abc4}(y, pg_t)=(fy,
g_t)=\exp\left\{it\left(\sum_{j=1}^as_jf_1 e_j
+\sum_{j=1}^l{q_j\over q}f_1 d_j +\sum_{j}k_jf_1 b_j
\right)\right\}.
\end{equation}
The required statement   follows from    (\ref{abc3})  and (\ref{abc4}).

Consider now the case when $K$
has infinite dimension.
Let  $N=\mathbb{R}^{\aleph_0*}\times\mathbb{Z}^{\mathfrak{m}*}$,
  ${\mathfrak{m}}\le\aleph_0$. Put
$M=\mathbb{R}^{\aleph_0}\times\mathbb{T}^{\mathfrak{m}}$.
The groups $N$ and $M$ are the character group of one another.

We need some properties of
nuclear and strongly  reflexive topological Abelian groups
(see \cite{Ba}).
We use them in the proof of Theorems 2 and 3.

The group
$N$ is   nuclear (\cite[(7.8), (7.10)]{Ba}). For such groups the Bochner  theorem
about one-to-one correspondence between the family of all
continuous positive definite functions
 on the group $N$ and  the family of all regular finite Borel
measures on its character group holds.  We note that the factor-group of
a nuclear group
with respect to a closed subgroup is also nuclear.

The groups  $M$ and $N$ are
 strongly  reflexive (\cite[(17.3)]{Ba}). These groups have, in particular, the following properties
 analogues  to that of locally compact Abelian groups.
The Pontryagin  duality theorem holds for such groups. Let $P$ be a closed subgroup of the group
    $M$.  For any $x\in M\backslash P$ there exists
 a character $y\in A(Y, P)$ such that $(x, y)\ne 1$. Any character of the subgroup $P$ can be extended to a character of the group $M$.
 Moreover, the natural homomorphisms
 $N/A(N, P)\mapsto P^*$ and $(M/P)^*\mapsto A(N, P)$  are topological isomorphisms.

\textbf{Theorem 2.} \textit{Let $X$ be a second countable  locally compact Abelian group such
that its connected component of zero  $c_X\cong\mathbb{R}^a\times K\times
\mathbb{T}^{\mathfrak{m}}$, where $a\geq 0$, $K$ be
a connected compact Abelian group containing no subgroups
topologically isomorphic to
 $\mathbb{T}$, ${\mathfrak{m}}\le \aleph_0$.
Assume that   $K$ has infinite dimension. Let $\alpha_j, \ \beta_j \in \mathrm{Aut}(X)$.
Let $\xi_j$, $j = 1, 2,\dots, n$, $n \ge 2,$ be
independent random variables with values in  $X$ and
 distributions  $\mu_j$ with non-vanishing characteristic
functions.
Assume that the linear forms   $L_1 = \alpha_1\xi_1 + \cdots +
\alpha_n\xi_n$  and $L_2 = \beta_1\xi_1 + \cdots + \beta_n\xi_n$
are independent.   Then there exist  a continuous monomorphism
$p:G\mapsto X$,  where   $G=
\mathbb{R}^{\aleph_0}\times\mathbb{T}^{\mathfrak{m}}$, and elements
 $x_j\in X$, \ $j=1, 2, \dots, n,$ such that all distributions
 $\mu_j*E_{x_j}$ concentrated on the subgroup
  $p(G)$. Furthermore, $p(G)$ is a characteristic subgroup.}

 \textbf{Proof.} We prove Theorem 2 in analogy with the proof of Theorem 1.
 As in the proof of Theorem 1 we can assume that $X=\mathbb{R}^a\times K\times
\mathbb{T}^{\mathfrak{m}}$, where     $K$ is
a connected compact Abelian group containing no subgroups
topologically isomorphic to
 $\mathbb{T}$, and  $Y=\mathbb{R}^a\times
D\times \mathbb{Z}^{\mathfrak{m}*}$, where  $D=K^*$
is a countable discrete torsion free Abelian
group containing no free factor-groups. Since the group
   $K$ has infinite dimension, the rank of the group
 $D$ is $\aleph_0$. Put
 $G=\mathbb{R}^{\aleph_0}\times \mathbb{T}^{\mathfrak{m}}$.
Then $H=G^*\cong\mathbb{R}^{\aleph_0*}\times
\mathbb{Z}^{\mathfrak{m}*}$.   We   assume that
$H=\mathbb{R}^{\aleph_0*}\times \mathbb{Z}^{\mathfrak{m}*}$.
Construct the natural embedding of the group    $Y$ to the group
$H$ (see e.g.
(\cite[\S5.9]{Fe})). For this purpose we choose in $D$ a  maximal
independent system of elements
   $d_1, \dots, d_l, \dots$ and, as we constructed in the proof of
   Theorem 1 the continuous homomorphism
  $f_0:{\mathbb{R}^a\times D\mapsto
\mathbb{R}^b}$, construct the  continuous homomorphism
$f_0:\mathbb{R}^a\times D\mapsto \mathbb{R}^{\aleph_0*}$.
Next extend it to the continuous monomorphism
 $f:Y\mapsto H$.  Since the group  $K$   has no subgroups topologically isomorphic to
  $\mathbb{T}$, we have
$\overline{f_0({\mathbb{R}^a\times D})}=\mathbb{R}^{\aleph_0*}$
(\cite[\S5.17]{Fe}). This obviously implies that $\overline{f(Y)}=H$.
Put $p=\tilde{f}$.  Since $G$ is a strongly reflexive group, it follows
from
   $\overline{f(Y)}=H$ that the homomorphism
$p:G\mapsto X$  is a monomorphism.

Consider the distributions   $\nu_j=\mu_j*\bar\mu_j$  and reason
as in the proof of Theorem 1.
We show that the restrictions of the characteristic functions
  $\hat\nu_j(y)$ to the subgroup $\mathbb{R}^a\times D$
are the characteristic functions of some Gaussian distributions
on the group $\mathbb{R}^a\times K$. Hence, the representations
 (\ref{new1}) and (\ref{new2}) hold,  where $Q_j$ are  infinite symmetric positive semidefinite matrices.
 Next, reason as    in the proof of Theorem 1 we come to the continuous positive definite functions
  $\bar h_j(w, k)$,
$w\in \mathbb{R}^{\aleph_0*}$, $k\in
\mathbb{Z}^{{\mathfrak{m}}*}$,  on the group $H$.
Since $H$ is a nuclear group,  we can correspond to each function  $\bar h_j(w, k)$
a  distribution  $\lambda_j$ on the group  $G$ such that
$\hat\lambda_j(w, k)=\bar h_j(w, k)$, $w\in
\mathbb{R}^{\aleph_0*}$, $k\in \mathbb{Z}^{{\mathfrak{m}}*}$.
We have $\hat\lambda_j(f(y))=\hat\nu_j(y)$   $y\in Y$. Taking into account that
 $G$ is a strongly reflexive group, it
implies that $\nu_j=p(\lambda_j)$.
The final part of the proof of Theorem 2 is the same as in Theorem 1.
  $\Box$

\textbf{Remark 2.}  Reasoning as in Remark 1 we see that in the
   case when the subgroup   $K$ has infinite dimension, the subgroup $p(G)$
   is also the arcwise connected component of zero of the group $X$.

We need the following general statement.

\textbf{ Proposition 1.} \textit{Let $X$ and $G$ be complete separable metric Abelian groups,
$\alpha_j, \ \beta_j \in \mathrm{Aut}(X),$ $j=1, 2, \dots, n,$
$n\ge 2,$ and  $p:G\mapsto X$ be a continuous monomorphism such that
 $\alpha_j(p(G))=\beta_j(p(G))=p(G)$, $j=1, 2, \dots n$. Let $\xi_j$
  be
independent random variables with values in the group $X$ and
 distributions  $\mu_j$. Assume that there exist elements $x_j\in X$
 such that all distributions $\mu_j*E_{x_j}$
are concentrated in the subgroup $p(G)$. Assume that the linear forms  $L_1 = \alpha_1\xi_1 + \cdots +\alpha_n\xi_n$
and $L_2 =
\beta_1\xi_1 + \cdots+ \beta_n\xi_n$ are independent.  Then
$\hat\xi_j=p^{-1}(\xi_j+x_j)$ are independent random variables
with values in the group
$G$, $\hat\alpha_j =p^{-1}\alpha_jp$,
$\hat\beta_j=p^{-1}\beta_jp$ are topological automorphisms of
the group $G$ and the linear forms  $\hat L_1  = \hat\alpha_1\hat\xi_1
+ \cdots
+\hat\alpha_n\hat\xi_n $ and $\hat L_2  =\hat\beta_1\hat\xi_1   + \cdots
+\hat\beta_n\hat\xi_n$ are independent}.

\textbf{Proof}. We use the following theorem by Suslin
\cite[\S 39, IV]{Ku}: Let $X_1$ be a complete separable metric space,
  $X_2$ be a metric space, $p: X_1\mapsto X_2$ be a continuous one-to-one
mapping. If $B$ is a Borel set in $X_1$, then $p(B)$ is a Borel
set in $X_2$ By the Suslin theorem  $\hat\xi_j$ are independent random
variables with values in the group $G$. It is obvious that the random variables
  $\hat\xi_j$  are independent. Since $p$ is
a continuous monomorphism and  $G$ is a complete separable
metric group, by the Suslin  theorem
 images     of Borel sets under the  mapping
   $p$ are also Borel. Hence $\hat\alpha_j$,
$\hat\beta_j$ are Borel automorphisms of the group    $G$, and hence
  $\hat\alpha_j$, $\hat\beta_j\in \mathrm{Aut}(G)$ (see e.g.
  \cite[\S 4.3.9]{Sr}). To proof the independence of the linear
  forms  $\hat L_1  $ and $\hat L_2 $
  it suffices to note that for any Borel subsets $A_1$ and $A_2$ of
  the group $G$ we have:
$\{\omega:\hat L_i(\omega)\in A_i\}= \{\omega:L_i(\omega)\in
p(A_i)-(\alpha_1x_1+\cdots+\alpha_nx_n)\}$,
 $i=1, 2$. $\Box$

 \textbf{ Remark 3.} Let $X=\mathbb{R}^a\times K\times
\mathbb{T}^{\mathfrak{m}}$,  where $a\geq 0$, $K$ be a connected compact
Abelian group containing no subgroups topologically isomorphic to
 $\mathbb{T}$ and having dimension      $l$,
${\mathfrak{m}}\le \aleph_0$.
Let $\xi_j, \ j = 1, 2,\dots, n, \ n \ge 2,$ be
independent random variables with values in  $X$ and
 distributions  $\mu_j$ with non-vanishing characteristic
functions.
Let $\alpha_j, \ \beta_j \in \mathrm{Aut}(X)$.
Assume that the linear forms   $L_1 = \alpha_1\xi_1 + \cdots +
\alpha_n\xi_n$  and $L_2 = \beta_1\xi_1 + \cdots + \beta_n\xi_n$
are independent. By Theorem 1 there exist
a continuous monomorphism
  $p:G\mapsto X$, where    $G=
\mathbb{R}^{b}\times\mathbb{T}^{\mathfrak{m}}$, $b=a+l$, and elements
 $x_j\in X$, \ $j=1, 2, \dots, n,$  such that all distributions
 $\mu_j*E_{x_j}, \  j = 1, 2,\dots, n,$ are concentrated on the subgroup
  $p(G)$, and  $p(G)$ is a characteristic subgroup.  It follows from Proposition 1 that
the study of
possible distributions of independent random variables
   $\xi_j$  is reduced to the study of
 possible distributions of independent random variables
$\hat\xi_j$ with values in  a group of the form
$\mathbb{R}^{b}\times\mathbb{T}^{\mathfrak{m}}$, where  $b\geq 0$,
 ${\mathfrak{m}}\le
\aleph_0$.

If the subgroup $K$ has infinite dimension,  applying Theorem
 2  instead of Theorem
  1, we come to the conclusion that
 the study of
possible distributions of independent random variables
   $\xi_j$  is reduced to the study of
 possible distributions of independent random variables
$\hat\xi_j$ with values in  a group of the form
$\mathbb{R}^{\aleph_0}\times\mathbb{T}^{\mathfrak{m}}$,
where  ${\mathfrak{m}}\le
\aleph_0$.

\bigskip

\noindent\textbf{3.  The Skitovich--Darmois theorem for
groups containing no subgroup topologically isomorphic to $\mathbb{T}^2$}

\bigskip

We will solve Problem 1 when  $n=2$  for second countable locally
compact Abelian groups  $X$ containing a subgroup topologically
isomorphic to $\mathbb{T}$
and containing no subgroups topologically isomorphic to
$\mathbb{T}^2$.  Taking into account Lemma 5, we can assume that  $X$
is a connected group. By Lemma    3  such group is topologically isomorphic to
the group of the form $\mathbb{R}^a\times K\times \mathbb{T}$, where $a\geq 0$, $K$ is a connected compact Abelian
group containing no subgroups topologically isomorphic
to   $\mathbb{T}$.
To avoid introducing
additional notation we assume that
$X=\mathbb{R}^a\times K\times \mathbb{T}$.  As has been noted in the
proof of
Theorem   1, $\mathbb{T}$ is a characteristic subgroup of the group
  $X$. Let  $\delta\in {\rm Aut}(X)$.
Since ${\rm Aut}(\mathbb{T})=\{\pm I\},$ we have either  $\delta|_\mathbb{T}=I$,
or $\delta|_\mathbb{T}=-I$.
Let $\xi_1, \xi_2$ be
independent random variables with values in  $X$ and
 distributions  $\mu_j$.
 Let $\alpha_j, \ \beta_j \in \mathrm{Aut}(X)$.  It is easily
 seen that the study
 of possible distributions
   $\mu_j$ provided that the linear forms
 $L_1=\alpha_1\xi_1+\alpha_2\xi_2$  and
$L_2=\beta_1\xi_1+\beta_2\xi_2$ are independent is reduced
to the case when
 $L_1$ and $L_2$ are of the form  $L_1=\xi_1+\xi_2$  and
$L_2=\xi_1+\delta\xi_2$, where  $\delta\in {\rm Aut}(X)$.
Then the Skitovich--Darmois functional equation (\ref{2})  becomes
\begin{equation}
\label{12} \hat\mu_1(u+v)\hat\mu_2(u+\varepsilon
v)=\hat\mu_1(u)\hat\mu_1(v)\hat\mu_2(u)\hat\mu_2(\varepsilon
v),\quad u,v \in Y,
\end{equation}
where  $\varepsilon=\tilde\delta$.

Let ${t}=(t_1,\dots, t_n, \dots)\in
\mathbb{R}^{\aleph_0}$ and  ${s}=(s_1,\dots, s_n, 0, \dots)\in
\mathbb{R}^{\aleph_0*}$.
Put $$\langle{t},{s}\rangle=
\mathop{\sum}\limits_{j=1}^\infty t_js_j, \quad
({t}, {s})=\exp\{i\langle{t}, {s}\rangle\}.$$
Let $\mu$ be a distribution on the group
   $\mathbb{R}^{\aleph_0}$.
The characteristic function of the distribution   $\mu$
is defined by the formula
 $$
\hat\mu({s})=\int\limits_{\mathbb{R}^{\aleph_0}}({t},
{s})d\mu({t}), \ {s}\in \mathbb{R}^{\aleph_0*}.$$  We remind that
a distribution  $\gamma$  on the group
$\mathbb{R}^{\aleph_0}$  is called Gaussian if its characteristic function
is represented in the form
$$\hat\gamma({s})=({t}, {s})\exp\{-\langle A {s}, {s}\rangle\},
\quad {s}\in \mathbb{R}^{\aleph_0*},
$$
where ${t}\in \mathbb{R}^{\aleph_0}$,
and $A=(\alpha_{ij})_{i, j=1}^\infty$ is  symmetric positive semidefinite matrix
   (see e.g. \cite[\S5.8]{Fe}).

We need some lemmas. It is convenient for us to formulate as a lemma
the following standard statement.

\textbf{Lemma 6.} {\it Let $X$ be a locally compact Abelian group,
$H$ be a closed  subgroup of   $Y$ and  $\mu \in {M}^1(X)$.
If $\hat \mu (y)=1$ for $y \in H$, then
the characteristic function $\hat \mu (y)$ is $H$-invariant and   $\sigma (\mu) \subset
A(X,H)$. }

{\textbf{Lemma 7}} (\cite{MiFe5}). \textit{ Let $Y$ be an Abelian group, $\varepsilon$ be an automorphism of the group   $Y$.
Assume that the functions $f_j(y)$ satisfy equation
$(\ref{12})$ and conditions   $f_1(0)=f_2(0)=1$.
Then each function $f_j(y)$ satisfies the equation}
\begin{equation}
\label{14} f_j(u+v)f_j(u-v)=f_j^2(u)f_j(v)f_j(-v),\quad u \in
(\varepsilon-I)Y, \ v\in Y.
\end{equation}

{\textbf{Lemma 8}} (\cite{Fe-Be-1986}). \textit{ Let $X$ be a second countable
locally compact Abelian group containing   no subgroups topologically isomorphic to
$\mathbb{T}^2.$ Let $\xi_1, \xi_2$ be
independent identically distributed random variables with values
in  $X$ and
 distribution  $\mu$. Assume that the  characteristic
 function $\hat\mu(y)$ does not vanish. If the linear forms
  $L_1=\xi_1+\xi_2$ and $L_2=\xi_1-\xi_2$ are independent, then  $\mu\in\Gamma(X)$.}

{\bf Lemma  9} (\cite{FeJFA}, see also  \cite[Lemma 12.4]{Fe5}).
{\it Let $Y$ be a topological  Abelian group,
let $\psi(y)$ be a continuous function  on $Y$ satisfying the
equation
\begin{equation}
\label{e6}\Delta_{2k}\Delta^2_{h} \psi(y) = 0, \quad h, k, y \in
Y,
\end{equation}
 and the conditions  $\psi(-y)=\psi(y), \ \psi(0)=0.$ Then
 the function $\psi(y)$ can be
represented in the form
\begin{equation}
\label{nnn1}\psi(y) = \varphi(y) + c_\alpha, \quad y \in y_\alpha +
\overline{Y^{(2)}},
\end{equation}
where $\varphi(y)$  is a continuous function on $Y$ satisfying
equation $(\ref{1})$, and $$Y = \bigcup_\alpha {(y_\alpha +
\overline{Y^{(2)}})}, \ y_0=0,$$ is a decomposition of the group
$Y$
 with respect to the subgroup
$\overline{Y^{(2)}}$.}

{\textbf{Lemma 10}}. \textit{ Let $X=\mathbb{R}^{\aleph_0}$,
$\alpha_j$, $\beta_j$, $j=1,2,\dots,n$, $n\geq 2$, be topological automorphisms of the group   $X$.
Let $\xi_j$ be
independent random variables with values in  $X$ and
 distributions  $\mu_j$.
Assume that the linear forms $L_1=\alpha_1\xi_1+\cdots+\alpha_n\xi_n$
and $L_2=\beta_1\xi_1+\cdots+\beta_n\xi_n$ are independent.
Then all $\mu_j$ are Gaussian distributions.}

\textbf{Proof.} Each element  ${s}\in \mathbb{R}^{\aleph_0*}$
defines a linear continuous functional  in the space
 $\mathbb{R}^{\aleph_0}$, and hence for any distribution
 $\mu$  on the group $\mathbb{R}^{\aleph_0}$
we can consider its image  $s(\mu)\in {M^1}(\mathbb{R}).$
We shall say that a distribution  $\gamma$ on the group
$\mathbb{R}^{\aleph_0}$ is weak Gaussian if
$s(\gamma)\in \Gamma(\mathbb{R})$ for all ${s}\in
\mathbb{R}^{\aleph_0*}$.
It is easy to see that the definitions of a Gaussian distribution and
a weak Gaussian distribution on the group
    $\mathbb{R}^{\aleph_0}$  are equivalent.    As has been noted in   \cite{My-2008},
the  Skitovich--Darmois theorem holds true for weak Gaussian distributions
on the group $\mathbb{R}^{\aleph_0}$.
$\Box$

{\textbf{Lemma 11}}. \textit{Let $X=\mathbb{R}^{\aleph_0}\times\mathbb{T}$.
Let $\xi_1, \xi_2$ be
independent identically distributed random variables with values in  $X$ and
 distribution  $\mu$. Assume that the  characteristic function $\hat\mu(y)$ does not vanish.
  If the linear forms
 $L_1=\xi_1+\xi_2$ and $L_2=\xi_1-\xi_2$ are independent, then  $\mu$
is a Gaussian distribution.}

\textbf{Proof.} The group $X$ is strongly  reflexive, and its
 character group is topologically isomorphic to the group  $\mathbb{R}^{\aleph_0*}\times\mathbb{Z}$.
and a Gaussian distribution on the group   $X$
is defined in the same way as   for the group $\mathbb{R}^{\aleph_0}$.
Since $X$ is strongly  reflexive, as easily seen, Lemma 1 holds for
 the group $X$.
In \cite[Lemma 9.11]{Fe1} Lemma 11 was proved  for second
countable locally compact Abelian groups
containing no more than one element of order 2. This proof is based on
Lemma 1 and  the properties of strongly  reflexive topological
Abelian groups listed before the formulation of Theorem 2. Hence, it
   is valid for the group
   $X$.
$\Box$

\textbf{Theorem 3.}  \textit{Let  $X=\mathbb{R}^a\times K\times
\mathbb{T}$, where $a\geq 0$, $K$ be a second countable connected compact Abelian group
containing no subgroups topologically isomorphic to
$\mathbb{T}$. Let $\delta \in {\rm Aut}(X)$.
Let $\xi_1, \xi_2$ be
independent random variables with values in  $X$ and
 distributions  $\mu_j$   with non-vanishing characteristic
functions.  Assume that the linear forms $L_1=\xi_1+\xi_2$
and $L_2=\xi_1+\delta\xi_2$
are independent. If $\delta|_\mathbb{T}=I$, then $\mu_j\in \Gamma(X)$.
If $\delta|_\mathbb{T}=-I$, then $\mu_j=\gamma_j*\pi_j$,
where $\gamma_j\in \Gamma(X)$, and $\pi_j$ are signed measures on the subgroup
$\mathbb{Z}(2)\subset \mathbb{T}$.}

\textbf{Proof.} There are two cases: either the dimension of the group
$K$ is finite or is infinite.

\textbf{1.}  Assume that the group  $K$  has a finite dimension  $l$.
Then the rank of the group   $D$ is also  $l$.  Taking into account Proposition 1 and
Remark 3, it suffices to prove Theorem 3
 for the groups of the form
$X=\mathbb{R}^b\times \mathbb{T}$, where  $b\ge 0$,
because in the notation of Proposition 1
$\delta|_\mathbb{T}=\hat\delta|_\mathbb{T}.$
We have $Y\cong\mathbb{R}^b\times \mathbb{Z}$.    To avoid introducing
additional notation we assume that
$Y=\mathbb{R}^b\times \mathbb{Z}$. Denote elements of the group   $X$
by $(t, z)$, $t\in \mathbb{R}^b$, $z \in \mathbb{T}$,
and elements of the group  $Y$ by $(s, n)$, $s\in\mathbb{R}^b$,
$n\in\mathbb{Z}$. It is easy to see that each automorphism
$\delta\in{\rm Aut}(X)$  is determined by the matrix
\begin{equation}
\label{28e1} \left(\begin{matrix}\alpha &v\\ 0&\pm 1\end{matrix}\right), \quad
\alpha\in{\rm Aut}(\mathbb{R}^b), \quad v\in \mathbb{R}^b,
\end{equation}
and the automorphisms   $\delta$ and  $\varepsilon=\tilde\delta$
act on the groups    $X$ and $Y$  respectively by the formulas
$$
\delta(t, z)=(\alpha t,e^{i\langle v, t\rangle}z^{\pm 1}), \ (t,
z)\in X, \quad \varepsilon(s, n)=(\alpha s+nv, \pm n), \ (s, n)\in
Y.
$$

It is obvious that   $\mathbb{R}^b$ is a characteristic subgroup of
the group   $Y$, because  $\mathbb{R}^b$ is a connected component of
zero of  $Y$.
Put $L={\rm Ker} (I-\varepsilon)|_{\mathbb{R}^b}$  and first verify
that the
proof of Theorem 3  is reduced to the case when  $L=\{0\}$.

By Lemma  1  the characteristic functions $\hat\mu_j(y)$ satisfy
equation  (\ref{12}).
Put  $\nu_j=\mu_j*\bar\mu_j, \  \ j=1, 2$. Then
$\hat\nu_j(y)=|\hat\mu_j(y)|^2> 0$ and the characteristic functions
  $\hat\nu_j(y)$ also satisfy equation
 (\ref{12}). Consider the restriction of equation
   (\ref{12})  for the functions
$\hat\nu_j(y)$ to the subgroup  $L$. We have
\begin{equation}
\label{t3.1}
\hat\nu_1(u+v)\hat\nu_2(u+v)=
\hat\nu_1(u)\hat\nu_1(v)\hat\nu_2(u)\hat\nu_2(v),
\quad u,v \in L.
\end{equation}
Set $h(y)=\hat\nu_1(y)\hat\nu_2(y)$. It follows from (\ref{t3.1})
that the function   $h(y)$ on the group  $L$ satisfies the equation
 $h(u+v)=h(u)h(v)$, i.e. $h(y)$ is a character of the group
$L$. This implies that the restrictions of the characteristic functions
   $\hat\nu_j(y)$ to $L$
are also characters of the subgroup   $L$. Since $\hat\nu_j(y)>0$, $
y\in Y$, we have   $\hat\nu_1(y)=\hat\nu_2(y)=1$ for $y\in L$.
Applying Lemma   6, we get $\sigma(\nu_j)\subset G=A(X,L)$.
Inasmuch as $L$ is the kernel of a continuous linear operator
in the space
 $\mathbb{R}^b$,  $L$ is a closed subspace in $\mathbb{R}^b$.
 It follows from this that $G=W\times
\mathbb{T}$, where  $W$ is a  closed subspace in
  $\mathbb{R}^b$. It is obvious that
  $\varepsilon(L)=L$. Hence
$\delta(G)=G$. It is clear that if $L\ne\{0\}$, then
$G$ is a proper subgroup of   $X$.

Consider a family of subgroups   $B$ of the group $X$ having
the properties:

$(i)$ $B=V\times\mathbb{T}$, where $V$ is a  closed subspace in     $\mathbb{R}^b$;

$(ii)$  $\delta(B)=B$;

$(iii)$ $\sigma(\nu_j)\subset B$, $j=1, 2$.

Let $N$ be an intersection of all subgroups of the group $X$
having properties     $(i)$--$(iii)$.
Obviously, the subgroup $N$  also possesses  properties
$(i)$--$(iii)$ and $N$ is the  smallest  subgroup having these properties.
 To avoid introducing
additional notation we assume that $N=\mathbb{R}^c\times\mathbb{T}$, $c\le b$, and  $N^*=\mathbb{R}^c\times\mathbb{Z}$.
Put  $\beta=\delta|_N$.
If ${\rm Ker} (I-\tilde\beta)|_{\mathbb{R}^c}\ne \{0\},$  then given above reasoning shows that the subgroup
$A(N, {\rm Ker} (I-\tilde\beta)|_{\mathbb{R}^c})$ possesses properties
$(i)$--$(iii)$ and is a proper subgroup of    $N$, contrary to the construction.
Since $\nu_j=\mu_j*\bar\mu_j$,  it follows from Lemma 4
that the distributions $\mu_j$ can be replaced by their
shifts $\mu'_j$ in such a manner that
$\sigma(\mu'_j)\subset N$. It follows from what has been said that we can assume from the beginning
 without loss of generality that  $L=\{0\}$. Let us note that the condition  $L=\{0\}$ means that in
 (\ref{28e1}) $(I-\alpha)\in{\rm
Aut}(\mathbb{R}^b)$.

{\text{\boldmath $1a$}}. Assume that $\delta|_\mathbb{T}=I$. It means that   the matrix which corresponds to the automorphism
  $\delta$  is of the form
$\left(\begin{matrix}\alpha &v\\ 0& 1\end{matrix}\right)$. Put $M={\rm Ker}
(I-\varepsilon)$   and consider the restriction of equation  (\ref{12})
for the functions $\hat\nu_j(y)$ to the subgroup  $M$.
We come to equation  (\ref{t3.1}), but when  $u, v \in M$.  Reasoning as above we conclude that
 $\sigma(\nu_j)\subset F=A(X,M)$. We have
 $F=A(X,M)=A(X,{\rm Ker}
(I-\varepsilon))=(I-\delta)(X)$. It follows from $X=\mathbb{R}^b\times
\mathbb{T}$, where $b\ge 0$, and $L=\{0\}$ that  $F=(I-\delta)(X)\cong\mathbb{R}^b$.
Since     $\varepsilon(M)=M$,  the restriction of the automorphism    $\delta$ to the subgroup $F$
is a topological automorphism of the group   $F$. Let $\eta_j$
be  independent random variables
 with values in  $F$ and
 distributions  $\nu_j$. Since the characteristic functions
$\hat\nu_j(y)$ satisfy equation  (\ref{12}),  the linear forms
 $\tilde L_1=\eta_1+\eta_2$ and $\tilde
L_2=\eta_1+\delta\eta_2$ by Lemma 1 are independent.  This implies by the Ghurye--Olkin theorem that
$\nu_j\in
\Gamma(F)$. Next, applying the Cram\'er theorem about decomposition of a Gaussian distribution the space
   $\mathbb{R}^b$ and Lemma  4, we get  $\mu_j\in
\Gamma(X)$.

{\text{\boldmath $1b$}}.  Assume that $\delta|_\mathbb{T}=-I$. It means that that the matrix which corresponds to the automorphism
   $\delta$  is of the form
  $\left(\begin{matrix}\alpha &v\\ 0&- 1\end{matrix}\right)$. Put
$H=(I-\varepsilon)Y$. Since  $L={\rm Ker}
(I-\varepsilon)|_{\mathbb{R}^b}=\{0\},$ we have
$H=Y^{(2)}=\mathbb{R}^b\times \mathbb{Z}^{(2)}\cong
\mathbb{R}^b\times \mathbb{Z}$. It follows from Lemma 7 that each characteristic function
  $\hat\mu_j(y)$ satisfies equation
 $(\ref{14})$ when $u, v \in Y^{(2)}$. Since   $H\cong
\mathbb{R}^b\times \mathbb{Z},$ we have $H^*\cong\mathbb{R}^b\times
\mathbb{T}$. Taking into account that the group  $H^*$ contains no subgroups topologically isomorphic to $\mathbb{T}^2$,
and applying Lemmas   1 and  8 we conclude that the restrictions of the characteristic functions
   $\hat\mu_j(y)$  to the subgroup
$H$ are the characteristic functions of Gaussian distributions.
Thus, we have  on the subgroup   $H$ the representation
$$\hat\mu_j(y)=m_j(y) \exp\{-\varphi_j(y)\}, \ j=1, 2,$$
where $m_j(y)$ are characters of the subgroup
 $H$, and $\varphi_j(y)$ are continuous nonnegative functions on     $H$  satisfying   equation
 $(\ref{1})$. Replacing if   necessary the distributions
 $\mu_j$ by  their shifts we can assume that
\begin{equation}
\label{15} \hat\mu_j(y)=\exp\{-\varphi_j(y)\}, \quad  y\in H,
\quad j=1, 2.
\end{equation}

Next the proof of Theorem 3 is quite similar to the proof of this theorem for the group
 $X=\mathbb{R}\times \mathbb{T}$ given in (\cite{MiFe5}). Put
$$ l_1(y)={\hat\mu_1(y)\over{|\hat\mu_1(y)|}}, \quad
l_2(y)={\hat\mu_2(y)\over{|\hat\mu_2(y)|}}, \quad y\in Y,
$$
and verify that  $l_j(y)$ are characters of the group  $Y$.
It follows from the equality $\hat\mu_j(y)=|\hat\mu_j(y)|$ for $y\in H$,
that $l_j(y)=1$ for  $y\in H$. The functions  $l_j(y)$  satisfy equation
(\ref{14}), which takes the form
\begin{equation}
\label{n1} l_j(u+v)l_j(u-v)=1,\quad u \in H, \ v\in Y.
\end{equation}
Substitute in  (\ref{n1}) $u=(s/2,0), \ v=(s/2,n)$. We obtain
$$ l_j(s,n)l_j(0,-n)=1, \quad s\in\mathbb{R}^b, \ n\in\mathbb{Z}. $$
Multiplying both sides of this equation by      $l_j(0,n)$ and taking into account that
  $|l_j(y)|=1, l_j(-y)=\overline{l_j(y)}$, we find
\begin{equation}
\label{n2} l_j(s,n)=l_j(0,n), \quad s\in\mathbb{R}^b, \
n\in\mathbb{Z}.
\end{equation}
Obviously, the functions    $l_j(y)$  satisfy equation
(\ref{12}). Taking into account    (\ref{n2}), we can write equation (\ref{12}) for the functions $l_j(y)$ in the form
\begin{equation}
\label{n3}
l_1(0,m+n)l_2(0,m-n)=l_1(0,m)l_1(0,n)l_2(0,m)l_2(0,-n),\quad m,n \in
\mathbb{Z}.
\end{equation}
We can get by induction from equation (\ref{n3})  that
$l_j(0,n)$ are characters of the group   $\mathbb{Z}$, and hence $l_j(y)$
are characters of the group   $Y$. Thus, there exist elements
$x_j \in X$ such that
\begin{equation}
\label{n5} l_j(y)=(x_j, y), \quad y \in Y, \quad j=1, 2.
\end{equation}

  Now find the representations for   $|\hat\mu_j(y)|$.
Put $\psi_j(y)=-\log|\hat\mu_j(y)|$. Then   (\ref{12}) implies that the functions
  $\psi_j(y)$  satisfy the equation
\begin{equation}
\label{m1} \psi_1(u+v)+\psi_2(u+\varepsilon v)=P(u)+Q(v), \quad u,v
\in Y,
\end{equation}
where $P(u)=\psi_1(u)+\psi_2(u), \ Q(v)=\psi_1(v)+\psi_2(\varepsilon
v)$.

As has been proved in \cite[Lemma 10.9]{Fe1},  (\ref{m1}) implies that each of the functions $\psi_j(y)$
satisfies the equation
\begin{equation}
\label{m4} \Delta_{(I-\varepsilon)k}\Delta^2_h \psi_j(y)=0,
\quad y, k, h\in Y.
\end{equation}
Since  $(I-\varepsilon)Y=Y^{(2)}$,  we conclude from  (\ref{m4}) that  each of the functions $ \psi_j(y)$ satisfies the equation
\begin{equation}
\label{m5} \Delta_{2k} \Delta_{h}^2 \psi_j(y)=0, \quad y, k, h \in
Y.
\end{equation}
 Since $H=Y^{(2)}=\overline{Y^{(2)}}),$  the decomposition of the group
$Y$ with respect to the subgroup  $\overline{Y^{(2)}}$  is of the form
$Y=H\cup((0,1)+H)$. It follows from  (\ref{15})  and Lemma 9 that there exist real constants
    $c_1, c_2$ such that
\begin{equation}
\label{nnnn} \psi_j(y) =  \begin{cases}
\varphi_j(y), & y \in H,
\\  \varphi_j(y) + c_j, & y\in (0,1)+H,
\end{cases}
\end{equation}
$j=1, 2$.

It follows from (\ref{12})  and (\ref{15})  that the equality
\begin{equation}
\label{30e1} \varphi_1(u+v)+\varphi_2(u+\varepsilon
v)=\varphi_1(u)+\varphi_1(v) +\varphi_2(u)+\varphi_2(\varepsilon
v)
\end{equation}
holds true for any   $u,v \in H$, and hence, for any
$u,v \in Y$.
  Substituting in  (\ref{m1}) $u, v\in (0,1)+H$ and taking into account
 (\ref{nnnn}) and (\ref{30e1}),  we find that
$c_1=-c_2$. Put $c_1=-c_2=-2\kappa$. It follows from  (\ref{nnnn}) that the functions
$|\hat\mu_j(y)|$  are of the form
$$ |\hat\mu_1(y)|=\exp\{-\varphi_1(y)+\kappa(1-(-1))^n\},\quad y=(s,n)\in Y,$$
$$ |\hat\mu_2(y)|=\exp\{-\varphi_2(y)-\kappa(1-(-1))^n\},\quad y=(s,n)\in Y.$$
Taking into account (\ref{n5}) we finally obtain
\begin{equation}
\label{m10} \hat\mu_1(y)=(x_1,
y)\exp\{-\varphi_1(y)+\kappa(1-(-1))^n\},\quad y=(s,n)\in Y,
\end{equation}
\begin{equation}
\label{m11} \hat\mu_2(y)=(x_2,
y)\exp\{-\varphi_2(y)-\kappa(1-(-1))^n\},\quad y=(s,n)\in Y.
\end{equation}

Consider the signed measures
$$\pi_1={1\over 2}(1+e^{2\kappa})E_{(0, 1)}+{1\over
2}(1-e^{2\kappa})E_{(0, -1)}, \quad \pi_2={1\over
2}(1+e^{-2\kappa})E_{(0, 1)}+{1\over 2}(1-e^{-2\kappa})E_{(0, -1)}$$
supported in   $\mathbb{Z}(2)\subset X$ (it is clear that one of $\pi_j$ is a distribution). Obviously,
the characteristic functions
 $\hat\pi_j(y)$  are of the form
\begin{equation}
\label{m1111} \hat\pi_1(y)=\exp\{\kappa(1-(-1)^n)\}, \quad
\hat\pi_2(y)=\exp\{-\kappa(1-(-1)^n)\}, \quad y=(s, n)\in Y.
\end{equation}
The statement of Theorem 3 follows from
(\ref{m10}), (\ref{m11})  and (\ref{m1111}). Note that
$\pi_1*\pi_2=E_{(0, 1)}$.

\textbf{2.}   Assume that the group  $K$  has infinite dimension. It follows from Proposition 1 and
Remark 3 that it suffices to prove Theorem 3
   for the group  of the form
    $X= \mathbb{R}^{\aleph_0}\times\mathbb{T}$.  Then  $Y\cong\mathbb{R}^{\aleph_0*}\times \mathbb{Z}$.    To avoid introducing
additional notation we assume that
$Y=\mathbb{R}^{\aleph_0*}\times \mathbb{Z}$.
    We follow the scheme of the proof of Theorem 3  in the case {\bf 1}. First note that  Lemma 1 holds true for the group $X$ because $X$ is strongly reflexive. Moreover,
Lemma 6 is also valid for the group $X$ because $Y$ is a nuclear group and $X$ is  strongly reflexive.
As in the case {\bf 1}  we consider the subgroup      $L={\rm Ker}
(I-\varepsilon)|_{\mathbb{R}^{\aleph_0*}}$  and reduce the proof of Theorem 3
to the case when
 $L=\{0\}$. In so doing we use the fact that any closed linear subspace in
   $\mathbb{R}^{\aleph_0}$ is  either finite-dimensional or  topologically isomorphic to
 $\mathbb{R}^{\aleph_0}$
(\cite{{BP}}),  and  that Lemma 6 holds true for the group $X$.

{\text{\boldmath $2a$}}.  Let $\delta|_\mathbb{T}=I$. Reason as in the case   {\text{\boldmath $1a$}}.   We use the fact that
Lemma 1 is valid for the group $\mathbb{R}^{\aleph_0}$, and that   Lemma 6 holds true for the group $X$.
Instead of the Ghurye--Olkin theorem we apply Lemma 10, and also use the fact that the Cram\'er theorem
about decomposition of a Gaussian distribution holds true for the space  $X=
\mathbb{R}^{\aleph_0}$.

{\text{\boldmath $2b$}}. Let $\delta|_\mathbb{T}=-I$.  Reason as in the case
{\text{\boldmath $1b$}}. We use the fact that
Lemma 1 is valid for the group $X$ and $Y$ is a nuclear group. Instead of Lemma 8 we apply Lemma    11. In all other respects
the proof repeats the proof in the case
    {\text{\boldmath $1b$}}.
 $\Box$

\newpage

Gennadiy Feldman

B. Verkin Institute for Low Temperature Physics and Engineering

of the National Academy of Sciences of Ukraine

Nauky Ave. 47

61103 Kharkiv, Ukraine

e-mail: feldman@ilt.kharkov.ua

\vskip 2 cm

Margaryta Myronyuk

B. Verkin Institute for Low Temperature Physics and Engineering

of the National Academy of Sciences of Ukraine

Nauky Ave. 47

61103 Kharkiv, Ukraine

e-mail: myronyuk@ilt.kharkov.ua

\vskip 1 cm

\noindent The corresponding author: G.M. Feldman, e-mail:
feldman@ilt.kharkov.ua


\newpage

\small

\begin{thebibliography}{99}

\bibitem{Ar}  D.L. Armacost. {\it The structure of locally compact abelian
groups}. Marcel Dekker, Inc., New York and Basel, 1981.

\bibitem{Ba} W. Banaszczyk. {\it Additive subgroups of topological vector
spaces}.  Lecture Notes in Mathematics, {\bf 1466},
Springer-Verlag, Berlin, 1991.

\bibitem{BP} Cz. Bessaga, A. Pelczynski. {\it
On a class of $B\sb 0$-spaces}.
Bull. Acad. Pol. Sci. Cl. {\bf III},   (1957), 375--377.

\bibitem{ChisGets}  G. P. Chistyakov, F. G\"otze. {\it Independence of linear forms with random coefficients}.  Probab. Theory Relat. Fields, \textbf{137},   no. 1-2, (2007), 1--24.

\bibitem{TVP1979}  G. M. Feldman. {\it Gaussian distributions on locally compact Abelian groups}.  Theory Probab. Appl., {\bf 23}, no. 3, (1979),   529--542.

 \bibitem{Fe-Be-1986}  G. M. Feldman. {\it Bernstein Gaussian distributions on
groups}.  Theory Probab. Appl., \textbf{31}, no. 1, (1986), 40--49.

\bibitem{SMZh}  G. M. Feldman. {\it Characterization of the Gaussian distribution on groups by the independence of linear statistics}.  Siberian Math. J., {\bf 31}, (1990),   336--345.

\bibitem{Fe}   G.M. Feldman. \textit{
Arithmetic of probability distributions and characterization
problems on Abelian
    groups }.  (AMS translation of mathematical monographs
\textbf{ 116}),  Providence, RI, 1993.

 \bibitem{Fe-SD-2003} G.M. Feldman. {\it A characterization of the Gaussian
distribution on Abelian groups}. Probab. Theory Relat. Fields,
\textbf{126}, (2003), 91--102.


\bibitem{Fe1} G.M. Feldman.  {\it  Functional equations and characterization problems on locally
compact Abelian groups}. EMS Tracts in Mathematics {\bf 5},
European Mathematical Society (EMS), Zurich, 2008.

\bibitem{Fe5} G.M. Feldman. {\it Characterization problems of mathematical statistics on locally compact Abelian groups}.  Naukova Dumka, Kiev, 2010.

\bibitem{FeJFA} G.M. Feldman. {\it The Heyde theorem for locally compact Abelian groups}. J.
of Funct. Analysis,  {\bf 258}, (2010),  3977--3987.

\bibitem{JOTP2015} G.M. Feldman. {\it On the Skitovich--Darmois theorem for the group of p-adic numbers}. J. of Theor. Probab., \textbf{28}, no. 2, (2015), 539--549.
    
    \bibitem{FG} G.M. Feldman, P. Graczyk. {\it The Skitovich-Darmois theorem for locally compact Abelian groups}. J. of the Australian Mathematical Society,  {\bf 88}, no 3, (2010), 339--352.

    \bibitem{MiFe5} G.M. Feldman, M.V. Myronyuk. {\it Independent linear forms on
connected Abelian groups}. Math. Nach.  \textbf{284}, no. 2--3, (2011), 255--265.

\bibitem{Fu1}  L. Fuchs.   {\it Infinite Abelian groups}. {\bf 1},
Academic Press, New York, San Francisco and London, 1970.

\bibitem{Hewitt-Ross}  E. Hewitt, and K.A. Ross. {\it Abstract
   Harmonic Analysis}. {\bf 1}, Springer-Verlag, Berlin,
    Gottingen, Heildelberg, 1963.
    	
\bibitem{Ibrahim} I. A. Ibragimov. {\it On the Skitovich--Darmois--Ramachandran Theorem}.
Theory Probab. Appl., {\bf 57}, no. 3, (2013), 368--374.

\bibitem{Kag-Lin-Rao} A. M. Kagan, Yu. V. Linnik, C.R. Rao.
{\it Characterization problems in mathematical statistics}.
  Wiley Series in Probability and Mathematical
Statistics, John Wiley $\&$ Sons, New York-London-Sydney, 1973.

\bibitem{Ku} K. Kuratowski. {\it Topology}, {\bf 1}, Academic Press,
New York and London, 1966.

\bibitem{Ma} I.P. Mazur. {\it Skitovich--Darmois Theorem for Discrete and Compact Totally Disconnected Abelian Groups}. Ukrainian Math. J.,
 {\bf 65}, no. 7, (2013), 1054--1070.

\bibitem{My-2008}  M.V. Myronyuk. {\it On the Skitovich-Darmous theorem and
Heyde theorem in a Banach space}. Ukrainian Math. J., {\bf 60}, no. 9, (2008),   1437--1447.


\bibitem{MiFe4}   M.V. Myronyuk, G.M. Feldman. {\it Independent linear statistics
on a two-dimensional torus}.
 Theory Probab. Appl., \textbf{52}, (2008),  78--92.

 \bibitem{NeueScho} D. Neuenschwander, R. Schott. {\it The Bernstein and
Skitovich-Darmois characterization theorems for Gaussian
distributions on groups, symmetric spaces, and quantum groups}.
Exposition. Math. \textbf{15}, no. 4,  (1997), 289--314.


\bibitem{Parthasarathy} K.R. Parthasarathy.  {\it Probability measures
on metric spaces}. Academic Press, New York and London, 1967.

\bibitem{Sr}  S.M Srivastava. {\it  Course on Borel sets}.
 Springer-Verlag, New York, 1998.





\end{thebibliography}
\end{document}